# RICHARD DEDEKIND Y LA ARQUITECTURA DEL CONTINUO ARITMÉTICO

Luis Giraldo González Ricardo
*Universidad de La Habana – UH – Cuba*

Carlos Sánchez Fernández
*Universidad de La Habana – UH – Cuba*

**Resumo**

É comum considerar que a tendência estruturalista matemática começa no século XX, em algum momento depois do trabalho sobre los fundamentos da teoria dos conjuntos e torna-se a sua propagação através do grupo Bourbaki. Neste artigo argumentamos que essa tendência estilística estava presente em Richard Dedekind (1831-1916) desde 1854 na sua dissertação de habilitação como professor na Universidade de Göttingen. O objetivo principal deste artigo é mostrar como se desenvolve o estilo estruturalista nos trabalhos de Dedekind e argumentar por qué consideramos que é importante para comprendeer a arquitetura do continuum aritmético.

**Palavras-chave:** Dedekind, estruturalismo matemático, número.

**[RICHARD DEDEKIND AND THE ARCHITECTURE OF THE ARITHMETICAL CONTINUUM]**

**Abstract**

It is usually considered that the structuralist tendency in mathematics began in the twentieth century, at some point after the works on set theory and obtained its spreading through the works made by the Bourbaki group. In the present paper we argument the presence of this stylistic inclination in Richard Dedekind (1831-1916) when he made this dissertation for the habilitation as privatdozent at the University of Göttingen in 1854. Our main objective is to show how evolved the structuralist style in Dedekind´s works, and to argument why we consider him significant to the architecture of arithmetical continuum.

## 1. Introducción

Varios trabajos relativamente recientes [Reck, 2003; Sieg & Schlimm, 2005; Yap, 2009] resaltan el papel significativo de la obra de Dedekind en la formación del estilo estructuralista, que dominó primero en la matemática germana y que más tarde fuera difundido universalmente por el grupo Bourbaki.

Los rasgos característicos del estructuralismo matemático a los que hacemos referencia, según [Reck-Price, 2000], pueden resumirse por:

- La matemática está principalmente ligada a la investigación de las estructuras.
- Esta investigación matemática conlleva una abstracción de la naturaleza de los objetos individuales
- Los objetos cumplen solamente lo que se puede establecer en términos de las relaciones básicas de la estructura.
- Para el "working mathematician" es importante el desarrollo y utilización de la teoría de conjuntos como esqueleto (*framework*) y elemento unificador de la investigación matemática.

En este trabajo nos hemos propuesto responder la cuestión siguiente: ¿Cómo evoluciona el estilo matemático de Dedekind hasta lograr atributos bien cercanos a los rasgos del estructuralismo matemático del siglo XX? Trataremos de argumentar, además, que desde su habilitación como *privatdozent* en 1854 ya mostraba esta inclinación por el estilo estructuralista y que no lo abandona en toda su obra, sino que por el contrario lo precisa en su búsqueda afanosa de dar una versión adecuada de la arquitectura del continuo aritmético.[1]

## 2. Los esbozos de la arquitectura aritmética

Las diferentes opiniones que sitúan a Dedekind como un exponente de las ideas estructuralistas no hacen insistencia en sus primeros trabajos. En nuestra opinión, las ideas expuestas en su habilitación se deben tomar en cuenta como significativas. En lo que sigue trataremos de argumentar por qué le damos tanto valor. Y nos parece importante resaltar que estas ideas no tienen todavía la influencia directa de Dirichlet o Riemann, con razón considerados sus maestros, pero en la época posterior al deceso del septuagenario Gauss (1777-1855), que fuera su guía científico oficial sin brindarle apoyo sistemático. El mismo año 1854 Riemann y Dedekind, ambos patrocinados por Gauss, presentan tesis de habilitación por la Universidad de Göttingen. Dedekind, que tenía solo 23 años de edad, tituló la disertación *Über die Einführung neuer Funktionen in der Mathematik.* El tema central era la introducción de nuevas operaciones en el dominio de la matemática. Dedekind destaca cómo en la matemática la necesidad de generalizar conceptos, no elimina aquellos que han sido generalizados, sino que deben mantenerse como casos particulares. Los ejemplos utilizados por Dedekind para ilustrar los planteamientos de su disertación provienen de la aritmética. Hagamos una síntesis de sus ideas principales.

---

[1] Notas biográficas que no incluimos en este artículo y el esbozo de algunas de las ideas que aparecen desarrolladas con profundidad aquí, pueden verse en el libro [Sánchez-González, 2013].

La repetición del proceso más simple de hallar el sucesor da origen a la adición, la repetición de la adición origina la multiplicación; y esta última a hallar potencias. Pero la incapacidad de realizar las operaciones inversas da origen a la *creación* de los nuevos dominios numéricos, con tal de poder efectuar la sustracción y la división.

La *creación* de los nuevos dominios numéricos conlleva a definir las operaciones directas sobre los números recién creados. Las operaciones aritméticas sobre un dominio más amplio son distintas a las definidas sobre los enteros positivos. Esta afirmación realizada por Dedekind muestra su interés por el rigor y la concepción de lo que diferencia las operaciones cuando se cambia el dominio de definición de las mismas. Por esta razón es que pone énfasis en que la operación de adición sobre los enteros (o los racionales) restringida a los números naturales tiene que operar de igual forma que la operación adición definida a partir de la repetición del proceso de hallar el sucesor.

El objetivo de Dedekind es establecer un procedimiento para extender las operaciones existentes a nuevos dominios numéricos más amplios, sin alterar su significado aritmético.

Las ideas expuestas por Dedekind muestran su interés en la comprensión del papel primordial que juegan las operaciones, en este caso aritméticas, en su relación con los elementos sobre los que operan. Dedekind expone en su disertación que los números irracionales aparecen por la imposibilidad de obtener de los números racionales las potencias con exponente racional no entero. Pero en este caso destaca que la dificultad para extender a través de su método las operaciones de suma y multiplicación a los números reales y después a los complejos, es mucho mayor:

> "*Estos pasos hacia adelante son tan inmensos que es difícil decidir cuál de los diferentes caminos que se presentan ante nosotros debemos seguir. Obviamente, si se quisiera aplicar a estas nuevas clases de números –los irracionales e imaginarios- las operaciones de la aritmética como ellas han sido desarrolladas hasta ahora, entonces se*

> *hacen necesarias extensiones repetidas de las definiciones anteriores – suma y multiplicación- y entonces (…) surgen las principales dificultades de la aritmética.[2]"* [Dedekind, 1854 en Ewald, 1996 pág. 759]

Aunque esta cita nos pudiera parecer pesimista, la visión de Dedekind era positiva sobre la obtención de una rigurosa definición de las operaciones aritméticas en los números reales y en los complejos:

> *"Sin embargo, debemos esperar que se logrará un edificio realmente sólido de la aritmética, si aplicamos persistentemente el principio fundamental: no permitir arbitrariedades, sino siempre guiarnos por las leyes descubiertas."* [Dedekind, 1854 en Ewald, 1996 pág. 759]

La disertación ofrecida por Dedekind como ejercicio de habilitación muestra su interés en la sistematización de los métodos para extender las operaciones matemáticas de cierto dominio a uno más amplio. Sin embargo, aún no hay trazas notables sobre la utilización de herramientas propias de los métodos teórico-conjuntistas posteriormente introducidas por Georg Cantor y el propio Dedekind.

La tesis de habilitación de Dedekind permite afirmar que en los inicios de la labor investigativa, a la temprana edad de 23 años, Dedekind tiene preferencia por el tratamiento conceptual de la matemática y en privilegiar las operaciones con relación a los elementos. Estas preferencias muestran su inclinación hacia el estructuralismo, aunque todavía no use el estilo teórico-conjuntista, ni defina concretamente ninguna estructura aritmética.

En lo siguiente veremos como la obra de Dedekind se mantiene fiel, en línea general con lo planteado en su disertación de habilitación y cómo las ideas "estructuralistas" van precisándose a medida que su cultura matemática se amplía.

## 3. La *nueva ola* del pensamiento matemático germano.

La formación matemática de Richard Dedekind ocurrió en la Universidad de Göttingen en el período 1850-1852, fecha en que defendió su doctorado bajo los auspicios de Gauss. A pesar de la enorme autoridad matemática que representa Gauss, ya éste había cumplido los 70 años de edad y no podía –ni le interesaba- dedicar su tiempo a la formación de matemáticos en Göttingen. Por lo que podemos afirmar que el nivel matemático de la Universidad de Göttingen era deficiente, y aún más, en todo el territorio de la antigua Confederación Germánica, solo la Universidad de Berlín tenía un nivel cercano, aunque inferior, a las instituciones universitarias de Francia que marchaban a la vanguardia, no solo de las ciencias exactas sino también técnicas. Sin duda, los temas centrales de la matemática de la época se manejaban fuera de Göttingen.

---

[2] Las traducciones a lo largo del artículo son responsabilidad de los autores.

Por el motivo antes mencionado el entrenamiento de Dedekind en temas fundamentales de la matemática de la época ocurrió después de obtener su habilitación. Tras la muerte de Gauss en 1855, su puesto fue suplido por Lejeune Dirichlet (1805-1859), quien enseguida se convirtió en maestro y amigo de Dedekind. A esto hay que sumarle la compañía de Bernhard Riemann (1826-1866), cinco años mayor que Dedekind, quien después de su llegada de una estancia larga en Königsberg y de la habilitación se convertiría también en colega y maestro. Subrayemos que esta influencia fue intensa, pero no extensa, ya que ambos, Dirichlet y Riemann, enfermaron y murieron pocos años después.

Dirichlet se formó, matemáticamente hablando, en la escuela francesa. Su novedoso tratamiento de las series trigonométricas de Fourier y la aplicación del análisis infinitesimal a la teoría de números le dieron un renombre dentro de la comunidad matemática europea. Dirichlet se distinguía de sus contemporáneos por preferir la determinación de las ideas matemáticas fundamentales por encima de los cálculos. Hermann Minkowski supo sintetizar brillantemente esta inclinación de Dirichlet:

> *"Superar los problemas matemáticos con un mínimo de cálculo y con un máximo de razonamiento perspicaz"* [Elstrodt, 2007]

Dedekind se convirtió en un asiduo alumno de las conferencias de Dirichlet, y en cartas a su familia mostró la influencia positiva que ejerció Dirichlet sobre él:

> *"(…) fue para mí un gran placer seguir sus conferencias profundas y agudas, seguí, en efecto, todas sus lecciones sobre la teoría de números, sobre el potencial, sobre la integral definida, sobre las ecuaciones en derivadas parciales, y él, tanto por sus enseñanzas como por las numerosas entrevistas personales que progresivamente devinieron más íntimas, hizo de mí un hombre nuevo."* [Dugac, 1976, pág. 21]

La acotación de *hombre nuevo* debemos verlo desde el punto de vista matemático, es claro que fue Dirichlet quien introdujo a Dedekind en las formas propias de la matemática más avanzada de la época. También se debe destacar que la asistencia de Dedekind a las lecciones de Dirichlet sobre teoría de números, y la posterior edición del texto del último sobre el tema le ofrecieron una motivación extra para realizar sus investigaciones sobre los números algebraicos.

Otra parte importante en la formación del pensamiento abstracto de Dedekind está vinculado al papel jugado por Riemann. Dedekind y Riemann, no obstante la diferencia de edad, hicieron amistad gracias al seminario de matemática y física dirigido por Wilhelm Weber y Moritz Stern en Göttingen, que ambos frecuentaban.

Dedekind también asistió a las clases impartidas por Riemann en Göttingen; entre los temas tratados estaban las integrales elípticas y la teoría de funciones de variable compleja. La visión que tenía Riemann de la matemática era bastante diferente a la estándar en la época; era, más bien, cercana a Dirichlet. Los métodos introducidos y desarrollados por Riemann eran más conceptuales que algorítmicos, como lo demuestran sus memorias póstumas –editadas y publicadas por Dedekind- *"Über die Darstelbarkeit einer funktion durch eine trigonometrische Reihe"* y *"Über die Hypothesen welche der Geometrie zu Grundeliegen"*. Para más detalles de la obra de Riemann, y sus relaciones con Dedekind, ver [Laugwitz, 1999].

El propio Dedekind afirmó:

> *"Además, tengo mucho contacto con mi excelente colega Riemann, quien es después, o incluso junto con Dirichlet, el más profundo de los matemáticos vivos, y así será reconocido pronto, cuando su modestia le permita publicar ciertas cosas que, sin embargo, serán solo comprensibles para unos pocos."* [Elstrodt, 2007]

La gestación de los métodos propios de Dedekind para hacer matemática, encontraron el ingrediente esencial en las nuevas ideas promovidas por Dirichlet y Riemann. La búsqueda de elementos conceptuales por encima de laboriosos cálculos se convirtió en ley en la obra de Dedekind.

Uno de los grandes logros de la obra de Dedekind lo constituye su nivel de abstracción y generalización, lo que permitió la extensión y aplicación de sus métodos a diversos problemas. Una de las premisas para la consecución de tal resultado es tener la habilidad de aislar el argumento matemático fundamental, del problema que se trate, la capacidad de sintetizar lo esencial del asunto que se estudia; tales competencias fueron desarrolladas en Dedekind gracias al intercambio con matemáticos renovadores como Dirichlet y Riemann, y no se mostró lerdo en usarlas en el campo que había abierto magistralmente su primer tutor Karl F. Gauss.

## 4. La primera elaboración de la arquitectura de los números algebraicos.

La expansión de los dominios de la aritmética a los números complejos mostró desde temprano su poder: Euler logró la demostración del teorema de Fermat para el caso *n=3* extrapolando las propiedades de **Z**, a una nueva clase de números: {x=p+q$\sqrt{-3}$: p y q enteros}. A principios del siglo XIX Gauss utilizó los números de la forma $a+bi$, con $a$ y $b$ enteros –los enteros gaussianos- para la prueba de la ley de reciprocidad bicuadrática.

Gauss, a diferencia de Euler, no extrapoló las leyes de la aritmética de los números enteros, sino que probó que en los enteros gaussianos se puede establecer teoremas análogos a los existentes en la aritmética tradicional, como es el caso del teorema fundamental de la aritmética.

Para la demostración de numerosos casos particulares del teorema de Fermat, el alemán Ernst Eduard Kummer (1810-1893) extendió las ideas de Gauss al caso de los enteros ciclotómicos. Con el estudio de este caso, diferente al de los enteros gaussianos, Kummer mostró que no siempre se puede extender la aritmética tradicional a los nuevos dominios aritméticos. La forma en que se solventó dicho problema fue ingeniosa: Kummer introdujo los *factores ideales* que permitieron recuperar la factorización única, y así extender las leyes usuales de la aritmética.

Estos trabajos de principios del siglo XIX dejaban claro que el estudio de las propiedades de los números enteros –lo que constituye el espacio de la aritmética tradicional- necesitaba la expansión a dominios más generales; la aritmética pura promovida por Gauss en sus *Disquisitiones Arithmeticae* encontró en ella dificultades insalvables. El mismo Gauss y después Kummer mostraron que era imprescindible cambiar la mentalidad para promover una nueva aritmética. Los trabajos de Kummer y Gauss no fueron casos aislados, también Dirichlet y F. Eisenstein utilizaron extensiones de los números enteros para probar el teorema de Fermat para el caso $n=5$, el primero; y la ley de reciprocidad cúbica el segundo.

¿Podían estas ideas renovadoras extenderse al estudio general de los enteros algebraicos? Este problema fue atacado desde diferentes perspectivas por un discípulo directo de Kummer, Leopold Kronecker (1823-1891), y por nuestro Dedekind. El primero siguió por caminos más trillados, aunque sin dejar de ser original[3], pero la vía radicalmente novedosa sería presentada por Dedekind con su estilo estructuralista.

La solución elaborada por Dedekind vio la luz en 1871, como el X suplemento de la segunda edición de *Vorlesungen über Zahlentheorie* de Lejeune Dirichlet.

Antes de ahondar en los métodos empleados por Dedekind hagamos algunos comentarios sobre la solución establecida por Kummer para la factorización de los enteros ciclotómicos, y así poder establecer los elementos novedosos de la obra de Dedekind.

La metodología seguida por Kummer fue profundamente algorítmica. La introducción de los factores ideales estaba muy ligada al dominio **Z**[α] –α es una raíz primitiva de la unidad- en el que se trabajara. Los números ideales primos se introducen de forma específica a cada caso, es decir, los factores ideales se obtienen a partir de las representaciones que se tienen de los mismos; en vez de utilizar características esenciales. Además, los métodos para hallar dichos factores ideales eran sumamente laboriosos. Otra dificultad en la obra de Kummer es que no aparece una definición explícita de lo que es un factor ideal, su existencia se verifica a través de los números que divide.

---

[3] Para más detalles en la solución de Kronecker al problema de la factorización única, ver [Edwards, 1990].

La primera edición de la teoría de Dedekind está, aún, basada en numerosos cálculos, (ver detalles en [Edwards, 1980]). La principal diferencia entre Kummer y Dedekind es que el último realiza una definición explícita del concepto de ideal, al lograr aislar las características esenciales de los factores ideales de Kummer. La definición de ideal, según Dedekind, está basada en ciertas reglas que deben satisfacer los elementos:

> *"Un sistema **J** de infinitos elementos de **D** será llamado ideal si satisface las dos condiciones siguientes:*
>
> *1.       La suma y la diferencia de dos elementos de **J**, es de nuevo un elemento de **J**.*
>
> *2.       Cada producto de un elemento de **J**, con un elemento de **D**, es un elemento de **J**."* [Dedekind, 1871, pág. 452]

Dedekind prueba, y no de forma sencilla, que a cada factor ideal de Kummer le corresponde un único ideal, y viceversa. De esa forma prueba que los resultados de Kummer se mantienen como caso particular de su teoría.

Pero no es solo esta la diferencia entre Kummer y Dedekind. Los pasos revolucionarios de Dedekind se sienten desde el inicio del X suplemento:

> *"Por un cuerpo* (Körper) *entiendo un sistema infinito de números reales o complejos, que es cerrado y completo en sí mismo, es decir, la adición, sustracción, multiplicación y división de cualesquiera dos de estos números siempre produce un número del mismo sistema"*
> [Dedekind, 1871 , pág. 424]

El aspecto más novedoso en la obra de Dedekind es la utilización, aunque de forma aún ingenua, de métodos conjuntistas. Para desarrollar su teoría Dedekind se ve necesitado de introducir nuevos conceptos, como el ya visto de cuerpo e ideal; pero también introduce la noción de Z-módulo y anillo[4]; y todos ellos son ciertos sistemas de números que cumplen propiedades específicas correspondientes.

Es necesario comentar que la noción de cuerpo no es completamente original de Dedekind, pues fue introducida –aunque de forma difusa- por Galois en sus estudios de la solución por radicales de las ecuaciones algebraicas. No se puede negar de forma categórica la influencia de las ideas de Galois, pues Dedekind era un conocedor profundo de la teoría de Galois; no por casualidad fue Richard Dedekind quien impartió el primer curso sobre teoría de Galois en las universidades de Alemania[5] [Dugac, 1976]. Sin embargo, en [Kiernan, 1971] se afirma que la fuente principal para el estudio de la teoría de cuerpos en Alemania tiene su origen en las *Disquisitiones Arithmeticae* de Gauss.

---

[4] El término anillo fue acuñado por Hilbert en el *Zahlbericht*.

Por otro lado, en toda la extensión del X Suplemento no se utiliza ninguna representación particular de los enteros algebraicos, lo que significa que Dedekind trabajó con ellos de forma independiente a las características concretas de cada dominio de enteros algebraicos. Lo anterior se traduce en que Dedekind se ha abstraído de las representaciones de los enteros algebraicos. Tal abstracción que realiza Dedekind de la expresión de entero algebraico **le permite centrarse en las operaciones entre los elementos de los sistemas en cuestión, más que en los elementos mismos**.

Es válido señalar que los resultados obtenidos por Dedekind, a pesar de su generalidad, se restringen al caso de subcuerpos de los números complejos, lo que es bastante claro, dado que no existen definiciones abstractas de los conceptos introducidos por Dedekind, nótese que la definición dada por Dedekind solo incluye los subcuerpos de **C**, en particular su interés se centra en el estudio de las extensiones finitas de **Q**. El nivel de generalidad lograda en la obra de Dedekind posibilita el tratamiento unificado de los enteros algebraicos de extensiones finitas del cuerpo de los números racionales. Eso sí, la notación empleada permite la extensión sencilla de sus resultados al caso abstracto.

Otro aspecto novedoso incorporado por la obra de Dedekind es la utilización del método axiomático fuera de la geometría, fruto de ello es la definición de *ideal* a partir de sus características fundamentales.

Por las razones antes expuestas podemos concluir que los métodos utilizados por Dedekind comienzan a separarse de los utilizados por sus predecesores, con la introducción de herramientas conjuntistas en su investigación matemática, aunque sin jugar todavía el rol principal dentro de la misma.

En 1871 todavía no se puede afirmar que la obra de Dedekind tenga todos los rasgos del estructuralismo matemático: las ideas teórico conjuntistas están circunscritas a subconjuntos de los números complejos, y sus métodos no han probado ser de validez fuera del contexto de la teoría de números. A pesar de lo antes dicho, la obra de Dedekind es completamente revolucionaria, y como tal fue recibida con frialdad por la comunidad matemática en general. Muchos de los matemáticos de la época no lograron comprender lo novedoso de los métodos utilizados por Dedekind. En carta a Lipschitz, Dedekind muestra la fría acogida de sus logros:

---

[5] Otros tempranos conocedores de la teoría de Galois en Alemania fueron E. E. Kummer y L. Kronecker [Bashmakova-Smirnova, 2000].

> "*Con excepción del profesor H. Weber de Königsberg, (…) quien ha expresado su interés en familiarizarse con la teoría* [de ideales. NA] *Ud. es el primero, no solo en mostrar interés en la materia, sino también de una forma tan práctica, que revive mis esperanzas de que mi trabajo no haya sido en vano*" [Dedekind, 1996, pág. 45]

## 5. La depuración de la idea de estructura.

La versión de 1871 no fue la única ofrecida por Dedekind de la teoría de ideales, en cada edición subsiguiente de las *Vorlesungen über Zahlentheorie* incluyó una nueva elaboración de su teoría. La tercera y cuarta edición vieron la luz en 1879 y 1894, respectivamente. Además de estas dos versiones, Dedekind publicó una versión en francés en *Bulletin des Sciences Mathématiques et Astronomiques* dirigido por Gastón Darboux, a solicitud de Rudolf Lipschitz en 1876.

La primera diferencia es que la versión francesa de la teoría de ideales contiene menos volumen de cálculo que la anterior, aunque es mayor en extensión. Este punto marca un mayor alejamiento de Dedekind de los métodos utilizados por Kummer.

La segunda redacción de la teoría de ideales contiene varios cambios respecto a la versión original de 1871. La multiplicación de ideales se introduce prácticamente al final de la versión de 1872, en la sección 6 de §163 [Dedekind, 1871, pág. 459], e incluso es innecesaria para la consecución del objetivo fundamental; sin embargo, para 1876 la multiplicación de ideales se define al inicio de la memoria. Enseguida define lo que se debe entender por ideal primo y establece como teorema fundamental que todo ideal es el producto de ideales primos, mientras que en la primera versión el teorema equivalente está dado en términos del mínimo común múltiplo de los ideales que lo dividen[6].

**La segunda presentación de la teoría de ideales está más acorde a los principios propuestos por Dedekind en su ejercicio de habilitación**. La introducción de los ideales está motivada por la imposibilidad de factorizar de forma única los enteros algebraicos contenidos en una extensión finita de **Q**.

Si vemos los resultados de 1872 a la luz de los principios establecidos en 1854 tenemos que la incapacidad de obtener factorizaciones de forma única sobre los enteros algebraicos provoca la introducción –o creación- de nuevos "números", pues los ideales se comportan de forma similar a los números enteros: se pueden sumar y restar, además se pueden multiplicar con las mismas propiedades que los números enteros, en especial la fundamental propiedad de factorización única. De esta forma las leyes aritméticas de **Z**, se trasladan "literalmente" al conjunto de los ideales.

Las líneas directrices establecidas en 1854 se convierten en el *leit motiv* de la obra

---

[6] También se puede ver [Edwards, 1980] y [Avigad, 2004].

posterior de Dedekind.

La tercera edición de las *Vorlesungen über Zahlentheorie*, de 1879 contiene, en esencia, la formulación de la teoría de ideales presentada al público francés tres años antes.

Las herramientas introducidas por Dedekind para traer calma sobre los rebeldes números algebraicos mostraron su verdadero poder cuando atacaron y resolvieron un problema, -aparentemente- distinto.

En 1882 apareció un artículo escrito a cuatro manos entre Richard Dedekind y Heinrich Weber –ambos editores de las obras de Riemann- con el título "*Theorie der algebraischen Funktionen einer Veränderlichen*". El objetivo de la obra era obtener las demostraciones de los teoremas obtenidos por Riemann sobre la teoría de las funciones algebraicas, entre ellos el teorema de Riemann-Roch, que es fundamental dentro del estudio de las funciones algebraicas.

> "*Los trabajos más abajo expuestos tienen la finalidad de fundamentar la teoría de las funciones algebraicas de una variable, una de las principales creaciones de Riemann, de una manera sencilla, rigurosa y totalmente general a la vez.*" [Dedekind, 1882, pp 238]

Muchas de las demostraciones de Riemann apelaban a argumentos geométricos y de continuidad –entre ellos el llamado "principio de Dirichlet"[7]-, que se presentaban, de cierto modo, *extraños*. Para Dedekind y Weber las demostraciones vinculadas a las funciones **algebraicas**, debían ser puramente algebraicas.

La idea empleada por Dedekind y Weber fue explotar las similitudes existentes entre las funciones algebraicas y los números algebraicos. De tal forma pudieron trasladar los conceptos introducidos en las *Vorlesungen über Zahlentheorie*, prácticamente palabra por palabra al campo de las funciones algebraicas, y los conceptos propios de la teoría tenían su paralelo en los números algebraicos.

La memoria de Dedekind y Weber no tuvo una acogida proporcional a la importancia de los resultados obtenidos. Otra vez los novedosos métodos *totalmente generales* preferidos por Dedekind –y también por Weber- chocaban contra la corriente principal de la labor matemática.

A pesar de la tibia acogida del trabajo conjunto de Dedekind y Weber, mostró el poder de las idea introducidas por Dedekind poco más de una década atrás. Las investigaciones del dueto Dedekind-Weber mostraron como los números algebraicos y las funciones algebraicas de una variable comparten características esenciales, lo que en lenguaje moderno sería: los números algebraicos y las funciones algebraicas comparten la misma estructura.

Dedekind y Weber supieron diferenciar entre las características particulares de cada campo de investigación para aislar las características generales que posibilitan un tratamiento unificado. Con este artículo se destaca la generalidad de los métodos de Dedekind, además que resalta su capacidad para abstraerse de la naturaleza de los elementos estudiados. Es, además, una muestra de la madurez del razonamiento "estructuralista" que Dedekind profesa.

---

[7] El "principio de Dirichlet" es también conocido como "principio variacional de Thomson-Dirichlet" ver p.e. Sánchez-Valdés (2004) pp. 75-77.

Para la cuarta edición la exposición de la teoría de ideales se vuelve incluso más abstracta. Entre los principales cambios está la incorporación de más resultados sobre los **módulos**, [Edwards, 1980]. Los resultados obtenidos sobre módulos le permiten ampliar la noción de ideal, y así introduce los ideales fraccionarios. La formulación de 1894 fue considerada por Emmy Nöther, -coeditora de las obras completas de Dedekind, y gestora en gran medida del nacimiento y desarrollo del álgebra abstracta a inicios del siglo XX- como la más elaborada, y su preferida.

Con la última presentación de la teoría de ideales Dedekind muestra la maestría que ha alcanzado en el manejo de herramientas conjuntistas por él introducidas para domar a los números algebraicos. El nivel de abstracción de las ideas de Dedekind lo comienzan a acercar más a la matemática elaborada en el siglo XX. Es decir, para 1894 Dedekind domina ampliamente las ideas abstractas que se volverán propias de la teoría algebraica de los números y seminales para la determinación de la arquitectura de los números; y es, además, capaz de aplicarlas con éxito para el desarrollo de una teoría aparentemente ajena al estudio de los números algebraicos. Es esta una de las premisas fundamentales para considerar a Dedekind un real predecesor de la corriente estructuralista matemática.

## 6. La estructura del continuo aritmético

El análisis infinitesimal encontró, finalmente, en el concepto de límite la fundamentación rigurosa de muchos de sus conceptos básicos, como el de continuidad en un punto y función derivada. La figura que representa el nuevo paradigma de rigor, iniciado con el siglo XIX, es Agoustin Louis Cauchy; aunque sin demeritar la labor de otros personajes ilustres.

Aunque Cauchy logró bastante en el campo de la fundamentación del análisis, también dejó varios problemas sin solucionar, que incluso lo llevaron a cometer sonados errores, como es el caso de la convergencia uniforme y la continuidad uniforme. Otro tema al que no dedicó especial atención fue el de los números reales.

En su archiconocido *Course d'Analyse* de 1821, Cauchy no define de forma explícita qué es un número real, ni establece sus características fundamentales. Sin embargo, tal vez motivado por la representación decimal de los números reales, afirma que todo número real es el límite de alguna sucesión de números racionales.

La afirmación de Cauchy es intuitivamente cierta, pero un estudio cuidadoso de su definición de límite abre nuevas interrogantes:

> *"Cuando los valores sucesivos de una variable se aproximan indefinidamente a un valor fijo, terminando por diferir del valor fijo en una cantidad tan pequeña como se quiera, a tal valor se le llama límite del resto"* [Grabiner, 2005 pág. 7]

Notemos como para la existencia del límite tiene que existir el valor fijo. Si tomamos la afirmación de Cauchy como definición tenemos un círculo vicioso, pues si consideramos una sucesión de números racionales que no converge a un número racional, entonces no tiene límite si no se ha definido primero el número irracional que es límite. Es decir, los números irracionales se definen a partir del concepto de límite, y este necesita de los números irracionales para tener consistencia lógica.

La necesidad de una definición rigurosa de los números reales fue notada desde las primera décadas del siglo XIX. Figuras como Bernhard Bolzano, Martin Ohm o William Hamilton intentaron, sin éxito, la definición rigurosa de los números irracionales a partir de los racionales.

No fue hasta la segunda mitad de la decimonovena centuria que varios matemáticos lograron construcciones precisas de los números reales. Matemáticos tan ilustres como Karl Weiertrass, Georg Cantor y Eduard Heine, así como uno menos conocido: el francés Charles Méray, además de Richard Dedekind, se relacionaron de forma satisfactoria con el problema. En [Dugac, 2003] y en [Sánchez-Valdés, 2004] se analiza el proceso de aritmetización del análisis, y se brinda más información sobre el tema que rebasa nuestros objetivos. Aquí nos interesa destacar el estilo estructuralista que Dedekind utiliza, similar al usado en sus investigaciones sobre los números algebraicos.

Dedekind concibió las ideas centrales de su construcción de los números reales en 1858, mientras preparaba un curso de Cálculo Diferencial, en su estreno como profesor del Politécnico de Zürich. Pero su publicación se demoró hasta 1872, cuando conoció de las publicaciones por parte de Cantor y Heine de sus respectivas construcciones de los números reales. Con su característica habilidad para penetrar en los fundamentos de la matemática Dedekind vio las debilidades presentes en las bases del análisis infinitesimal.

> *"(...)por primera vez me vi obligado a explicar los elementos del Cálculo Diferencial y sentí más agudamente que nunca antes la falta de una fundamentación realmente científica de la aritmética."* [Dedekind, 1872 en Dedekind, 1963 pág. 1]

El problema encontrado por Dedekind fue la utilización de argumentos geométricos para la demostración de que toda sucesión monótona y acotada tiene límite. Dedekind considera que la demostración de este teorema debe hacerse solo con herramientas aritméticas, y así evitar los argumentos "extraños" provenientes de la geometría.

La construcción de Dedekind se basa en la noción de cortadura, que es un aporte original de su intelecto y se sustenta en la estructura de orden en los números racionales:

> *"Dado un número a, el conjunto de los números racionales ($\boldsymbol{Q}$) se divide en dos clases $A_1$ y $A_2$. En la primera de ellas están todos los números $a_1$ tales que $a_1 < a$, y en la segunda todos los $a_2$ tales que $a_2 > a$."* [Dedekind, 1872 en Dedekind, 1963 pág. 6]

La idea fundamental es que todo número racional define una cortadura, pero hay cortaduras en el conjunto de los números racionales que no están definidas por números racionales. El otro aporte de Dedekind en su monografía, es la caracterización de la continuidad a partir de la idea de cortadura:

> *"Si todos los puntos de la recta se dividen en dos clases, tales que todo punto de la primera clase está a la izquierda de todo punto de la segunda clase, entonces existe un único punto que produce esta división de todos los puntos en dos clases, este punto corta la línea recta en dos partes"* [Dedekind, 1872 en Dedekind, 1963 pág. 11]

Si una cortadura no está definida por un número racional, entonces se *crea* un número irracional, que queda completamente determinado por la cortadura. Conceptualmente, según la solución propuesta por Dedekind, un número irracional no es un punto, sino que está determinado por un par de conjuntos infinitos de números racionales. Lo anterior muestra, nuevamente, la completa aceptación por parte de Dedekind del infinito actual.

La construcción de Dedekind concluye con la demostración de que el conjunto de las cortaduras cumple la caracterización de la continuidad, lograda por él utilizando la noción de cortadura, es decir, el conjunto de los números reales es un dominio continuo. Además define la suma para los números recién creados, y demuestra que toda sucesión monótona y acotada de números reales es convergente.

### ¿Se mantiene fiel Dedekind a las ideas adelantadas por él en su tesis de habilitación?

El problema que encuentra Dedekind, para buscar una construcción de los números reales está ligado a la imposibilidad de realizar la operación de pasar al límite en una sucesión monótona y acotada de números racionales, la que no siempre tendrá un número racional a la que se aproxime indefinidamente. El problema anterior Dedekind lo reduce a la imposibilidad de realizar la operación de corte en el conjunto de los números racionales.

El siguiente paso es la definición de las operaciones aritméticas para los nuevos números de forma que se mantengan inalteradas para los números racionales. En la monografía de 1872 Dedekind solo estudia el caso de la suma, y expone que la extensión del resto de las operaciones se hace de forma sencilla. Una elaboración exhaustiva de la teoría del número real según Dedekind se puede encontrar en el texto de Edmund Landau *Foundations of Analysis,* [Landau, 1966].

En *Stetigkeit und irrationale Zahlen* Dedekind se acoge a los principios elaborados en su tesis de habilitación 18 años atrás[8], además tiene varios puntos de contacto con la teoría elaborada para domesticar los números algebraicos. Primeramente, la construcción de los números reales no reside en ninguna representación particular, ni en elaboraciones calculistas. La esencia del método de Dedekind reside en la utilización en la característica de **Q**, de ser un conjunto totalmente ordenado y denso dondequiera. Este último detalle, junto a la no utilización de representaciones particulares permite la extrapolación del método de Dedekind al completamiento de cualquier conjunto denso y totalmente ordenado, prácticamente palabra por palabra.

Además de lo expuesto en el párrafo anterior, Dedekind mismo unificó el tratamiento por él dado a la creación de los ideales y los números irracionales; así lo muestra en una nota al pie de la memoria *Sur la Théorie des Nombres Entiers Algébriques* [Dedekind, 1996]. Dedekind se refiere a la necesidad de definir de forma clara y precisa la multiplicación de ideales, "*pues ellos no existen*" en el dominio de los enteros algebraicos:

> "*(…) la necesidad de tales demandas* –de establecer una definición de la multiplicación de ideales- *las que siempre deben ser impuestas con la introducción o creación de nuevos elementos aritméticos, se vuelve evidente cuando se compara con la introducción de números reales irracionales (…)*" [Dedekind, 1876 en Dedekind, 1996 pág. 57-58 ]

Estas ideas ya están presentes en el método matemático de Dedekind, al menos desde la consecución de su habilitación. No solo son puntos de contacto lo que existe entre la disertación de 1854 y la construcción de los números reales. En la más antigua Dedekind considera que la necesidad de introducir los números irracionales es la imposibilidad de hallar potencias de exponente racional de los números racionales, idea que es completamente abandonada en 1872. La explicación se encuentra en la nota al pie en

---

[8] Así lo destaca en la primera sección de [Dedekind, 1872].

[Dedekind, 1996] ya citada anteriormente:

> "(…) *se debe exigir que todos los números irracionales se engendren simultáneamente por una definición común, y no sucesivamente como raíces de ecuaciones, logaritmos, etc..*"

Luego de la discusión de los trabajos de Dedekind referentes a los números algebraicos y a los números irracionales, resalta la uniformidad de los métodos de Dedekind, para tratar problemas muy diferentes, aunque *en esencia* tratan ambos de la arquitectura del continuo aritmético.

### 7. Fundamentos de la Arquitectura de los Números.

Poco tiempo después de publicar *Stetigkeit und irrationale Zahlen* Dedekind comienza las labores para la creación de una monografía dedicada a elaborar una teoría que permitiese una definición rigurosa de los números naturales. La creación de la citada obra concluyó en 1888 con la publicación de *Was sind und was sollen die Zahlen?*.

Al mismo inicio del prefacio Dedekind plantea que "*en la ciencia nada debe aceptarse sin demostración*", lo que puede decirse que es el *leit motiv* de la memoria, a tal punto que logra una demostración del principio de inducción completa.

En el prefacio Dedekind adelanta la respuesta a la pregunta propuesta en el título -¿Qué son y para qué sirven los números?-:

> "*Considero el concepto de número completamente independiente de las nociones e intuiciones de espacio y tiempo, porque lo considero como un inmediato resultado de las leyes del pensamiento. Mi respuesta a los problemas propuestos en el título de esta monografía es, brevemente: los números son una creación libre de la mente humana; ellos sirven como medio para aprehender más fácilmente las diferencias entre las cosas. Es solo a través del proceso lógico de construir la ciencia de los números, y de esa forma adquirir el dominio aritmético continuo, en que estamos preparados adecuadamente para investigar nuestras nociones de espacio y tiempo, poniéndolos en relación con este dominio numérico creado en nuestra mente.*" [Dedekind, 1888 en Dedekind, 1963 pág. 32]

La cita muestra a las claras la posición de Dedekind, los números constituyen una herramienta para la investigación del mundo real; y son, además, creados por la mente humana. Por otra parte, el concepto de número no se origina en la noción de tiempo y espacio, como plantean Kant y sus seguidores, sino que

> *"Considero a la aritmética como una necesaria, o al menos natural, consecuencia del más simple acto aritmético, contar."*
> [Dedekind, 1872 en Dedekind, 1963 pág. 4]

La base sobre la que Dedekind construye su teoría del número natural es la teoría de conjuntos ("ingenua"), que elaboró de forma autónoma, aunque sin despreciar la influencia de Cantor al interrelacionarse ambos. Al desarrollo de las ideas básicas de la teoría de conjuntos dedica Dedekind la primera sección de *Was sind und was sollen die Zahlen?*.

La noción de conjunto que utiliza Dedekind es la típica de la infancia de la teoría, donde un conjunto se considera completamente determinado si cierto elemento pertenece o no al conjunto. Esta noción primitiva llevó poco tiempo después al desarrollo de las antinomias descubiertas por Cantor, Russell, y otros matemáticos de la época.

No es solo el desarrollo de los principios de la teoría de conjuntos, y las operaciones básicas entre conjuntos los temas novedosos introducidos por Dedekind; también podemos hallar la primera definición abstracta de función[9] –en la terminología de Dedekind: transformación (*abbildung*)-. El trabajo de Dedekind se centra en un tipo fundamental de transformaciones: las similares (*ähnlich*) que en lenguaje moderno serían las funciones inyectivas. Aunque en la definición dada por Dedekind se trate de funciones inyectivas la utilización que les da es propia de las funciones biyectivas, es decir, además de la correspondencia 1-1, considera que su imagen es todo el conjunto de llegada.

Las herramientas básicas utilizadas por Dedekind son las nociones de conjunto y función; pero la herramienta especializada es el concepto de *cadena* (kette). Un subconjunto K, de cierto universo A, es una cadena, respecto a la transformación φ, si φ(K) es un subconjunto de K. Nótese que esta definición no dice más que K sea estable por la transformación φ. Nuevamente vemos cierta similitud con el estudio de los números algebraicos, donde los ideales son conjuntos estables por la diferencia y el producto por un elemento externo; así como los cuerpos son estables por las cuatro operaciones aritméticas.

La noción de transformación similar le permite definir de forma rigurosa cuando un conjunto es infinito:

> *"Un sistema S se dice infinito si es semejante a una parte propia de sí mismo (...)"* [Dedekind, 1888 en Dedekind, 1963 pág. 63]

---

[9] Para un estudio más profundo del desarrollo del concepto de función ver [Sánchez-Valdés, 2007 ]

Una vez más el intelecto de Dedekind pone en práctica el "segundo principio de Dirichlet" para dar con la característica esencial de un conjunto infinito. En las secciones anteriores hemos visto como Dedekind usa de forma continuada conjuntos infinitos en sus investigaciones, lo mismo para la construcción de los números reales que para recobrar la factorización única en los enteros algebraicos.

Entre los conjuntos infinitos Dedekind distingue a cierta clase de conjunto infinito: los conjuntos simplemente infinitos. La naturaleza de los conjuntos simplemente infinitos no depende solamente de los elementos que los constituyen, sino también de la existencia de cierta transformación similar φ y un elemento especial, que se denota por 1. Luego, un conjunto simplemente infinito **N** está totalmente determinado por la existencia de una transformación φ de **N** y un elemento 1 que cumplan las siguientes condiciones:

α) $\varphi(N) \subset N$

β) $N = 1_0$

δ) $1 \notin \varphi(N)$     (*)

γ) φ es similar

Donde $1_0$ lo constituye la menor cadena que contiene a 1.

De la definición anterior se deben destacar dos puntos:

• El alto nivel de rigor lógico en la memoria de Dedekind: las características que definen a los conjuntos simplemente infinitos no vienen numeradas, sino se identifican con letras griegas. La razón es bien simple, ¡todavía no se han construido los números naturales!

• La definición de conjunto simplemente infinito, aunque no se diga de forma explícita, está dada por medio de axiomas. De nuevo las técnicas utilizadas por Dedekind lo suficientemente generales que se adaptan a diversas situaciones. Muestra la utilidad del método axiomático en diferentes áreas de la matemática, aunque en esta etapa del desarrollo matemático no se puede hablar de método axiomático formal como aparece después en Hilbert.

Una vez definidos los conjuntos simplemente infinitos Dedekind da su definición de los números naturales:

> *"Si en la consideración de un conjunto simplemente infinito **N** ordenado por una transformación φ negamos completamente el carácter especial de los elementos; simplemente reteniendo lo que los distingue entre si y tomando en cuenta solo las relaciones de uno con otro a través de la transformación ordenadora φ, entonces estos elementos son* *llamados* **números naturales** *o* **números ordinales** *o simplemente* **números** *(…)Con referencia a esta liberación de los elementos de cualquier otro contenido (abstracción) se justifica decir que los números son una creación libre de la mente humana. (…)"* [Dedekind, 1888, en Dedekind, 1963 pág. 68]

Entonces los axiomas dados para los conjuntos simplemente infinitos se convierten en axiomas para los números naturales. A la función φ se le llama ordenadora, y a la imagen de un elemento *n*, φ(*n*) se le llama sucesor de *n*. En el entorno de los conjuntos simplemente infinitos Dedekind demuestra el *Teorema de Inducción Completa*, el que unido a la transformación ordenadora se convierte en herramienta fundamental para construir la aritmética de los números naturales.

La definición de los números naturales descansa sobre la definición de la noción de conjunto simplemente infinito, ¿pero son todos iguales? De ser negativa la respuesta se pudiera obtener distintos "sistemas de números naturales". La respuesta de Dedekind es una demostración de "categoricidad", esto es: entre dos conjuntos simplemente infinitos existe una biyección que conserva la función ordenadora φ. Otra forma de decirlo es que todos los conjuntos simplemente infinitos son iguales salvo isomorfismo. Lo que en las palabras propias de la escuela de Hilbert, el único modelo que satisface las condiciones (*) es el conjunto de los números naturales[10].

El resto de la memoria se dedica al desarrollo de la aritmética usual de los números naturales, utilizando de forma magistral las propiedades de las funciones similares y el teorema de la inducción completa. En esta parte Dedekind muestra la seriedad con que pensó y elaboró su disertación de habilitación, pues sigue los pasos establecidos en 1854 para definir y establecer las propiedades fundamentales de la suma, la multiplicación y la potenciación con exponente natural.

En 1890, en carta dirigida al profesor de Hamburgo, Hans Keferstein aparece el siguiente comentario de Dedekind:

> *"¿Cuáles son las propiedades básicas, independientes entre sí, de esta serie N, es decir, aquellas propiedades que no pueden deducirse unas de otras pero de las cuales se siguen todas las demás? Y ¿de qué manera hay que despojar a estas propiedades de su carácter específicamente aritmético, de manera que queden subordinadas a conceptos más generales y a actividades del entendimiento tales que sin ellas no es posible en absoluto el pensamiento, pero con ellas viene dado el fundamento para*

---

[10] Para un análisis detallado del estudio de la categoricidad por Hilbert, y la influencia de Dedekind ver [Ortiz-Valencia, 2010].

> *la seguridad y completitud de las demostraciones, así como para la construcción de definiciones libres de contradicción?*" [Heijenoort, 1967 pág. 99-100]

La cita anterior muestra, una vez más, y de forma totalmente explícita la forma de comprender y resolver los problemas por parte de Dedekind: búsqueda de la idea fundamental del problema, abstracción y generalización; que son, además, elementos característicos de la actitud estructuralista en la matemática.

*Was sind und was sollen die Zahlen?* es, sin duda alguna, la obra de mayor abstracción redactada por Dedekind. Donde muestra todo su talento para extraer la esencia de los objetos estudiados, y sobre ello construir la teoría necesaria. Esta obra de 1888 es fruto de la madurez intelectual de Dedekind, en ella se funden las técnicas explotadas una y otra vez, así como nuevos métodos, como lo es la utilización de funciones para el estudio de las características de los números naturales.

En este punto de su carrera Dedekind se acerca mucho más a las líneas fundamentales de los matemáticos del siglo XX que comúnmente son considerados exponentes del estilo estructuralista.

## 8. Último aporte a la Arquitectura Aritmética[11]

En 1894 Dedekind se retiró de las labores docentes que realizaba en el Politécnico de su ciudad natal, pero no se retiró de la investigación matemática: en dos tardíos artículos, de 1897 y 1900, introduce la noción de retículo (*lattice*), que después sería redescubierto por Garret Birkhoff.

Los principios para el estudio y desarrollo de la teoría de los retículos se encuentran en fecha tan temprana como 1871, cuando define módulo en el X suplemento de las *Vorlesungen über Zahlentheorie*. Un módulo no es más que un subconjunto de los números reales o complejos estable para la adición y la sustracción, además, aplicando ideas propias de la teoría de conjuntos establece las relaciones entre divisor y múltiplo, máximo común divisor (mcd) y mínimo común múltiplo (mcm); lo que resulta natural de su propósito de extender la aritmética de los números enteros a los ideales.

En esta primera formulación de la teoría Dedekind no denota de forma particular las operaciones mencionadas en el párrafo anterior, pero con su exquisitez a la hora de exponer sus ideas, en la segunda versión de 1877 introduce símbolos específicos que le permiten esclarecer la naturaleza de los conceptos mencionados.

---

[11] En esta parte hemos utilizado fundamentalmente información brindada en el reciente artículo [Schlimm, 2011].

La notación más apropiada le permite notar cierta dualidad entre el mcd y el mcm. Tal fenómeno se conoce hoy en día como "leyes modulares" dentro de la teoría de retículos:

$(M_1)$ mcd[mcm(a, b), mcm(a, c)]=mcm[a, mcd(b, mcm(a, c))]
$(M_2)$ mcm[mcd(a, b), mcd(a,c)]=mcd[a, mcm(b, mcd(a, c))]
cualesquiera sean los números *a*, *b* y *c*.

En el proceso de edición de la obra de Dedekind, Emmy Nöther encuentra un manuscrito, no publicado, donde ya aparece un análisis de los fenómenos descritos arriba. La publicación de las ideas referentes a este tema verá la luz en los artículos:

**1897**: "*Über Zerlegungen von Zahlen durch ihre grössten gemeinsamen Teiler*"
**1900**: "*Über die von drei Moduln erzeugte Dualgruppe*"

Cuando generaliza las propiedades aritméticas del *mcd* y *mcm* utiliza dos operaciones generales denotadas por $\pm$ y plantea una definición abstracta de "*Dualgruppe*":

> *Si dos operaciones + y – sobre dos elementos arbitrarios A, B de un conjunto (finito o infinito) G generan dos elementos A±B del mismo conjunto G que satisfacen las condiciones:*
> *(1)*        *A+B=B+A  y A-B=B-A*
> *(2)*        *(A+B)+C=A+(B+C)  y (A-B)-C=A-(B-C)*
> *(3)*        *A+(A-B)=A y A-(A+B)=A*
> *entonces, sin interesar la naturaleza de estos elementos, G es llamado un Dualgruppe con respecto a las operaciones ±.* [Schlimm, 2011]

Tanto la notación como la terminología utilizada por Dedekind pueden confundir, pues el término *grupo dual* no es utilizado hoy en día con el mismo significado, como tampoco se *suman* o *restan* elementos del retículo. Expresado de la forma más usual ahora sería:

Un retículo es una estructura algebraica definida sobre un conjunto G por dos operaciones internas binarias $\Delta$ y $\nabla$ que cumplen los axiomas de conmutatividad (1), asociatividad (2) y absorción (3).

Los ejemplos estudiados por Dedekind muestran que en la estructura abstracta de Dualgruppe las propiedades (M₁) y (M₂) del *mcd* y el *mcm* son equivalentes a la única relación:

(M) [AΔ(B∇C)] ∇ (BΔC)=[A∇ (BΔC)] Δ (B∇C)
que hoy se conoce como "ley modular".

Dedekind demuestra también que la ley anterior es independiente de las 3 leyes que definen la estructura de Dualgruppe. El resultado mencionado permite a Dedekind introducir una nueva subestructura. Un Dualgruppe que verifica la propiedad (M) Dedekind lo llama "Dualgruppe de tipo modular".

Mostrando sus habilidades para trabajar con entes abstractos Dedekind construye otro tipo de estructura a través de las dos propiedades duales siguientes:

(D₁)  (A∇B) Δ (A∇C)=A∇ (BΔC)
(D₂)  (AΔB) ∇ (AΔC)=AΔ (B∇C)

que son llamadas "leyes distributivas". Para ellas prueba que son independientes a los axiomas que determinan el Dualgruppe y así se obtiene una nueva subestructura: un Dualgruppe que verifica ambas leyes distributivas Dedekind lo llama "Dualgruppe de tipo ideal".

Para estas dos nuevas subestructuras actualmente se usa el nombre de "retículo modular" y "retículo distributivo" respectivamente. Algunos autores en honor a su introductor le llaman "retículo de Dedekind", por ejemplo, en la literatura especializada rusa es común llamarlos "estructuras de Dedekind" (así aparece como "Dedekindova estructura" en la "Enciclopedia de Matemática", T.2 pág. 63, edición de 1979. Quizás valga aclarar que los nombres dados por Dedekind no son arbitrarios. Provienen de los dos ejemplos fundamentales que ha desarrollado en su obra anterior: los módulos de números enteros y los ideales de números reales algebraicos.

Inmediatamente Dedekind muestra que todo Dualgruppe de tipo ideal es de tipo modular y da ejemplo de Dualgruppe de tipo modular que no es de tipo ideal. Con esto prueba que los dos conceptos no coinciden, es decir que las leyes distributivas y las modulares no son equivalentes.

Emmy Nöther en los comentarios a este artículo poco conocido de 1897 dice:

"*es sobre todo significativo como una investigación axiomática adelantada a su tiempo*" [Schlimm, 2011]

y añade que es un ejemplo de construcción lógica a través de los modelos particulares apropiados para el sistema de axiomas en cuestión.

Poco tiempo después de publicar este artículo de 1897 Dedekind conoce la existencia de una monografía dedicada a un tema de lógica, pero donde se presenta una axiomática semejante a la introducida por él. La monografía se debe a Ernst Schröder, y lleva por título "*Vorlesungen über die Algebra der Logik*", en tres tomos. Los propósitos de Schröder y Dedekind son diferentes, pero muy similares desde el punto de vista estructural

Cuándo Dedekind conoce de la existencia de esta obra ya se había publicado los 3 volúmenes del "*Vorlesungen über die Algebra der Logik*" de Schröder y también el artículo de Dedekind de 1897. Con su espíritu analítico y sistemático se enfrasca en una lectura profunda hasta que comprende la semejanza entre el "cálculo lógico" de Schröder y sus "Dualgruppen". Decide publicar un nuevo artículo donde explícitamente plantea la conexión entre ambas teorías axiomáticas. Así surge su artículo de 1900 "*Über die von drei Moduln erzeugte Dualgruppe*". En más de 30 páginas Dedekind pretende mostrar el centro de su teoría y lo ilustra con un estudio de caso, una disección del retículo modular con 3 generadores. No solo llega a que su estructura está constituida exactamente por 28 módulos diferentes, sino que también determina la estructura del retículo distributivo generado por 3 ideales y llega a encontrar la cantidad reducida de elementos en este otro caso. Además prueba teoremas fundamentales sobre ambas estructuras.

Estas tardías obras de Dedekind, en plena madurez, muestran su capacidad de explotar al máximo el "segundo principio de Dirichlet" para encontrar la esencia del problema que investiga. Su razonamiento potente le permite encontrar las más profundas analogías y así aislar las características fundamentales que le permitan atacar varios problemas con los mismos recursos.

## 9. Comentarios finales

El interés primordial de Dedekind en toda su obra, más que la construcción de un edificio, es el establecimiento de sus cimientos y la distribución de las columnas que conforman su estructura. En síntesis, Dedekind en su laboreo con los números se comporta como un arquitecto, no como un simple albañil. Tal vez la mejor forma de comprender la profundidad del pensamiento de Dedekind, sea conocer la incomprensión que sufrió por parte de la mayoría de sus contemporáneos, y la admiración de la mayoría de los matemáticos posteriores.

En nuestra opinión la obra de Dedekind, aunque demorara en ser comprendida, sirvió de base para establecer, no solo una arquitectura del continuo aritmético, sino también la dirección principal en la fundamentación de la matemática en el siglo XX. Los principios esgrimidos por Dedekind son propios de la tendencia matemática hacia la formalización y el estudio de las relaciones funcionales entre las diferentes estructuras como interés primordial. En definitiva, la obra de Richard Dedekind sirvió de modelo para la posterior difusión del estilo estructuralista matemático, aunque no muchos lo reconocieran en sus publicaciones.

**Bibliografia**

**Luis Giraldo González Ricardo**
Departamento de Matemática.
Facultad de Matemática y Computación.
Universidad de La Habana. Cuba

**E-mail:** luis.gonzalez@matcom.uh.cu

**Carlos Sánchez Fernández**
Departamento de Matemática.
Facultad de Matemática y Computación.
Universidad de La Habana. Cuba

**E-mail:** csanchez@matcom.uh.cu